# Verification of the Proth-Gilbraith conjecture up to $\pi(10^{14})$ and methods for anticipating and minimizing verification results beyond that point

Simon Plouffe
October 24, 2025

Abstract

A calculation was performed to verify Proth-Gilbraith's conjecture for all prime numbers up to 10^14. The previous calculation was performed by Andrew Odlyzko in 1993 up to 10^13.

This involves calculating the differences between consecutive primes in absolute value and starting over. The conjecture states that all lines except the first begin with 1. To prove it, it suffices to find a line beginning with 1 and followed only by 0 and 2.

2, 3, 5, 7, 11, 13, 17, 19, 23, 29, 31, …
1, 2, 2, 4, 2, 4, 2, 4, 6, 2, …
1, 0, 2, 2, 2, 2, 2, 2, 4, …
1, 2, 0, 0, 0, 0, 0, 2, …
1, 2, 0, 0, 0, 0, 2, …
1, 2, 0, 0, 0, 2, …
1, 2, 0, 0, 2, …

The conjecture was verified by A. Odlyzko starting from the line composed of all prime numbers less than or equal to 10^13. A question suggested by Jean-Paul Delahaye was whether in 2025 we could do better and verify the conjecture for the sequence of prime numbers up to 10^14. Initially, I was able to go up to 2 ×10^13. I was then able to go up to 10^14. The calculation was completed on October 7, 2025, at 2:25 p.m. The article contains various information and statistics on the calculation and on the methods that can be used to anticipate the number of lines before arriving at the one containing a 1 followed by 0 and 2.

Preamble

A similar calculation was performed by Jean-François Colonna[*] at the same time, and he reached the same conclusion before me on October 5, 2025. The techniques used were C and Cshell, as well as a prime number generator based on the Miller-Rabin test. The conclusion of his calculation is the same: $G(\pi(10^{14})) = 693$. [15]

*CMAP (Centre de Mathématiques Appliquées), École polytechnique, France

## Introduction

Estimating the work to be done

There are 3,204,941,750,802 prime numbers up to 10^14 (3,204 billion), and we need to calculate the differences between these numbers up to the final line. A test with blocks of 10 million consecutive prime numbers taken from different locations shows that the average number of difference operations before reaching the final line is around 400. It is unreasonable to perform the calculation on all of these prime numbers in a single block. A good idea is to cut them into blocks of 10 million prime numbers and add an overlap for each block so that the calculations can be done independently. Knowing that the maximum number of operations per block is a few hundred, taking an overlap of several thousand values (3,000 to 4,000), we are sure to cover the entire range between 2 and 10^14. In this way, the calculations can be performed independently and/or in parallel on a machine. With this information, we can estimate that the total number of operations is around $1.2 \times 10^{15}$.

An initial attempt using Maple software showed that it was far too slow, requiring 271 days of continuous calculation on two computers. Following this failure, a program was created using Unix scripts and a series of commands native to this language. The result was almost as slow as Maple (241 days of calculation).

The alternative that was found is to use Python's numpy library (3.12) and the <diff> program, which allows the operation to be performed at very high speed. In addition, the size of the integers involved in the calculation can be controlled. Here, 16-bit numbers (max 32768) were used, which saves a considerable amount of time.

Equipped with this special program on two fast computers:

Intel I9-9900KF 8 cores/16 threads at 5 GHz and 32 gigabytes of memory
Intel I7-6700K 4 cores/8 threads at 4 GHz and 32 gigabytes of memory.

The difference operation on large arrays of numbers is performed at a speed of 2 billion calculations per second, which is 24 times faster than Maple or Unix scripts. In the end, the entire calculation took only 11.5 days.

The primes were generated with the *primesieve* program, which is the fastest known program for generating prime numbers up to $2^{64} - 1$. To save space, the calculation is performed on a just-in-time basis (explained below), with each file generated being processed immediately and only the result of the number of calculation steps to reach the final line being retained. Storing 3204 billion prime numbers requires 25 terabytes. Some calculations were performed by Maple or with a script to periodically check the validity of the calculation.

## Calculation plan

As suggested by Jean-Paul Delahaye, the work can be divided into slices, ensuring the validity of the operation by overlapping the lists of prime numbers. First, the width of a slice is chosen: here, 10 million primes. The size of the overlap is also chosen, between 3,000 and 4,000 prime numbers. By calculating the Proth-Gilbraith triangle (successive operations of absolute differences) of two successive slices, we reach the end of the calculation when the remaining line contains only 0s or 2s: only the first slice begins with 1. If the final slices, once the calculation is complete, overlap with enough values, then the calculation is validated. Empirically: the number of steps that Odlyzko obtained for his calculation was 635, which is the value of the G function: therefore, for $G(\pi(10^{13})) = 635$. This means that we should expect a few hundred steps at best. The maximum value encountered for $G(\pi(10^{14}))$ was 693. So, given that the overlap consists of more than 3,000 prime numbers, we can be sure that there are no gaps as long as the value of G does not exceed 1,000 or 2,000.

If the slices are 10 million primes to which we add 3000 or 4000 primes for the overlap, we therefore have 320495 slices of primes to process. A dynamic program has been developed to process the difference table on the fly, resulting in a trace of each slice once completed. In order, we have for each slice 10 million primes + overlap.

1- Creation of the list of primes.
2- *Awk* program to perform the difference between 2 primes on the fly or in a just-in-time flow, the program compares 2 successive lines as the first file is filled. The file is first generated and sent directly to the flow by the awk command but preceded by the pipe $< | >$.
3- Launch of the Python program with numpy.diff ( np.diff for short).
4- We keep track of the calculation to validate the number of steps required for each slice (320495 files).

Nomenclature:

G(range of prime numbers) = number of steps to get to a line containing only 0s and 2s, except the first slice which starts with 1.

G' (range of prime numbers) = number of steps to get to a line containing only 0 and 2. The calculation is performed on the largest known value M. Tables of record prime gaps between primes already exist. In this way the value of 693 was already detected on September 24 as the very possible final value of $G(\pi(10^{14}))$ because the value of M at an already known location is 766. Several sources confirm this maximum M after the prime number: 19581334192423 and which was discovered in 1989.

M(range of prime numbers) = maximum obtained following the difference operation carried out once. The value of M is the maximum between all the primes in the slice considered. Here the maximum encountered was 766 over the entire range between 2 and $10^{14}$.

R(range of prime numbers) = ratio between G(slice)/M(slice).

For example, G(slice # 25127) = 635 is the value found in 1993 by Odlyzko and which corresponds to the gap after the prime number 7177162511713. The value of 693 corresponds to the gap of 766 after the prime number 19581334192423. We have G(slice 66233) = 693 which is the maximum value over the entire range: $G(\pi(10^{14}))$.

## Results and statistics

We quickly realize that there is a direct link between the value of G and that of M, the R value is on average 0.830 and very few R values exceed 1: 99.2 % of the values are ≤1.

The maximum value found of 693 has as value of M which is 766, that is the maximum of the range after making the first difference between the primes. The ratio R is therefore 693/766 = .9046. This value was recently validated by the complete calculation of JF Colonna and the author of this article.

Another consequence of this average is that we can predict the value of G according to the value of M. Indeed, in table [7] we find the gap 906 after the number 218209405436543 found by T. Nicely in 1996. The value of G is 773. Here the ratio is 0.8532. We can also consult the list of known record gaps on Wikipedia (Prime gap) and other sources.

Here is a histogram of the R values calculated on 320,495 values of G and M. On the abscissa the value of R and on the ordinate the number of cases with this value. We quickly see that the vast majority of values are ≤1, more than 99.2%.

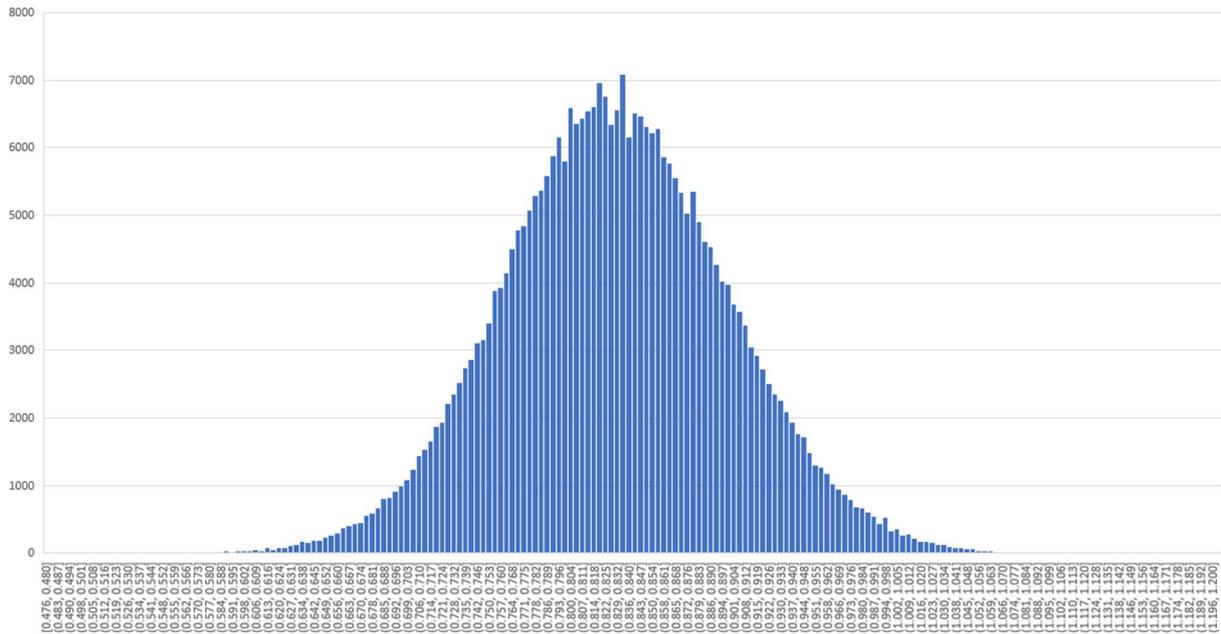

Graph of G values, we recognize the 2 known maximum points, the first at 635 as calculated by Odlyzko in 1993 (later rechecked) and the value of 693 detected now. The average G value is 375 on this sample.

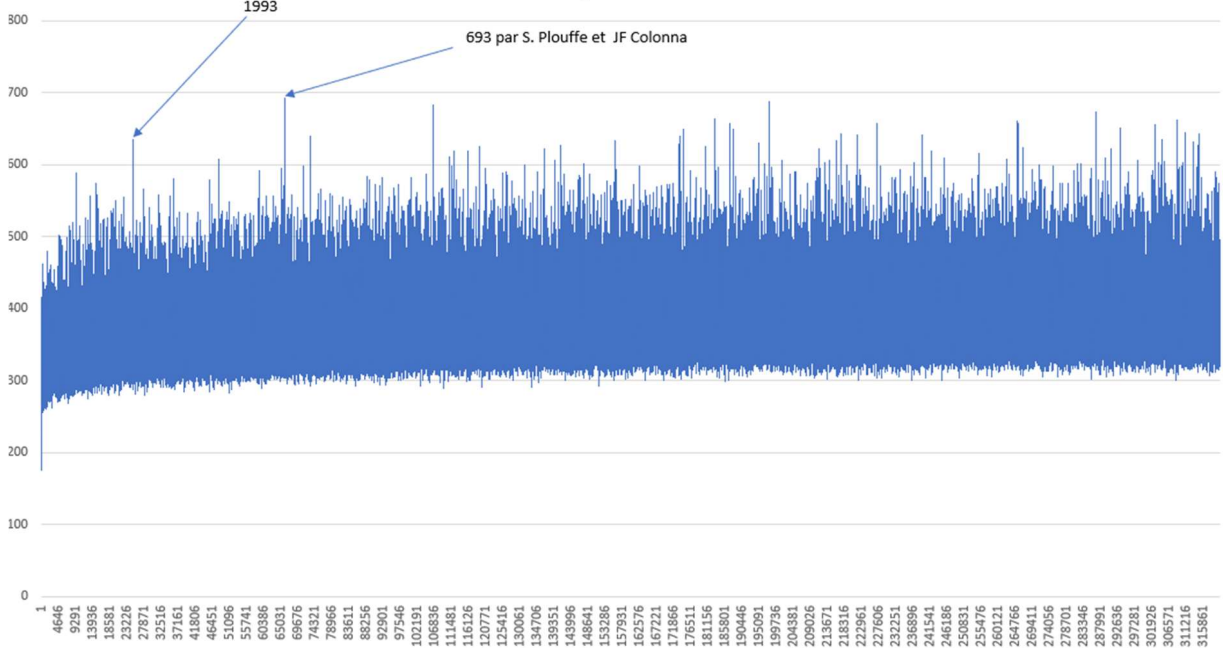

It is quite obvious that the value of G for a given slice or for the whole interval up to $\pi(10^{14})$ depends directly on the highest value M of this interval. Indeed, over the range of all primes up to $10^{14}$, the known maximum is 766 which is the largest known difference (Prime gap). Consequently, the calculation of G around this value actually gives the maximum value of G for the whole interval. So, knowing these values M: list of record prime gaps, it is quite easy to directly calculate the value G'(interval) which will probably be the final value for the whole interval. This is for a simple reason: For a given value of M the quantity G is on average (calculated on 320495 values) 0.830 and the maximum value exceeds 1 in only 0.8 % of cases. When the value exceeds 1 it is at most 1.20 so if we have the maximum known value M on a prime number there is a good chance that the calculation of G'(interval) gives the final value for G. The list in [13] was passed in full to calculate all the values of G'( $\pi(10^k)$), for k = 12 to 20. It was also noted that the known maximum values of M in table [13] correspond exactly to the values found by Odlyzko in 1993 [3].

Here is the table of values of G'( $\pi(10^k)$).

| K | Value of M | Value of G' | Ratio | Prime |
|---|---|---|---|---|
| 15 | 880 | 788 | 0.8954 | 277900416000927 |
| 16 | 1132 | 970 | 0.8560 | 1693182318746371 |
| 17 | 1184 | 1199 | 1,012 | 43841547845541059 |
| 18 | 1358 | 1274 | 0.9381 | 523255220614645319 |
| 19 | 1446 | 1433 | 0.9910 | 9656919634106230133 |
| 20 | 1402 | 1462 | 1,042 | 10103695526434940251 |

We could push the calculation further but to do this: we must have a list of gaps between the primes that are as small as possible, the list that appears on several sources allows us to go up to G'( $\pi(10^{20})$), after this mark it is more speculative. We have known gaps up to primes of 87 digits but we do not know if the prime number found is the smallest. Consequently, we can only give an indication of the minimum of G'( $\pi(10^m)$) where m > 20.